\documentclass[reqno, 12pt]{amsart}
\usepackage[letterpaper,hmargin=1in,vmargin=1in]{geometry}
\usepackage{amsmath, amssymb, amsthm, verbatim, bbm, url} 
\usepackage{graphicx}

\usepackage[OT2, T1]{fontenc}
\usepackage[russian, english]{babel}

\newcommand \N {\mathbb{N}}
\newcommand \R {\mathbb{R}}
\newcommand \C {\mathbb{C}}
\newcommand \Z {\mathbb{Z}}

\newcommand \Oh {\mathcal{O}}

\newcommand \D {\partial}
\newcommand \eps {\varepsilon}

\newcommand \Def {\stackrel{\textrm{def}}=}
\newcommand\comp{{\mathrm{comp}}}

\DeclareMathOperator \re {Re}
\DeclareMathOperator \im {Im}

\DeclareMathOperator \res {res}

\DeclareMathOperator \Hess {Hess}

\DeclareMathOperator \Tr {Tr}

\DeclareMathOperator \vol {vol}

\DeclareMathOperator \grad {grad}
\DeclareMathOperator \Div {div}
\DeclareMathOperator \spec {spec}
\DeclareMathOperator \loc {loc}

\theoremstyle{definition}

\numberwithin{equation}{section}
\numberwithin{lem}{section}
\numberwithin{Defn}{section}

\title
{Inverse problems in spectral geometry}

\author[Kiril Datchev]
{Kiril Datchev}
\author[Hamid Hezari]
{Hamid Hezari}
\address{Mathematics Department, Massachusetts Institute of Technology, Cambridge, MA
02139.}
\email{datchev@math.mit.edu}
\email{hezari@math.mit.edu}
\thanks{The first author is partially supported by a National Science Foundation postdoctoral fellowship, and the second author is partially supported by the National Science Foundation under
grant DMS-0969745. The authors are grateful for the hospitality of the Mathematical Sciences Research Institute, where part of this research was carried out.}
\date{August 29, 2011}

\begin{document}

\maketitle

\section{Introduction}

In 1966, Marc Kac in his famous paper \cite{Kac}  raised the following question:
Let $\Omega \subset \mathbb{R}^2$ be a bounded domain and let $$ 0\le \lambda_0 < \lambda_1 \leq \lambda_2 \dots $$ be the eigenvalues of the nonnegative
Euclidean Laplacian $\Delta_\Omega$ with either Dirichlet or Neumann boundary conditions. Is $\Omega$  determined up to isometries from the sequence $\lambda_0, \lambda_1 \dots$? We can ask the same question about bounded domains in $\mathbb R^n$, and below we will discuss other generalizations as well. Physically, one motivation for this problem is
identifying distant physical objects, such as stars or atoms, from the light or sound they emit. These inverse spectral problems, as some engineers in \cite{R, RWShN, PWR, RWP} have recently proposed, may also have some interesting applications in shape-matching, copyright and medical shape analysis.

The only domains in $\mathbb R^n$ which are known to be spectrally distinguishable from all other domains are balls. It is not even known whether or not ellipses are spectrally rigid, i.e. whether or not any continuous family of domains containing an ellipse and having the same spectrum as that ellipse is necessarily trivial. We can go further and ask the same question about a compact Riemannian manifold $(M,g)$ (with or without boundary); can we determine $(M,g)$ up to isometries from the spectrum of the Laplace-Beltrami operator $\Delta_g$? Or in general, what can we hear from the spectrum? For example, can we hear the area (volume in higher dimensions or in the case of Riemannian manifolds) or the perimeter of the domain? For the sake of brevity we only mention the historical background for the case of domains.  In 1910, Lorentz gave a series of physics
lectures in G\"ottingen, and he
conjectured that the asymptotics of the counting
function of the eigenvalues are given by: $$ N(\lambda)=\sharp\{ \lambda_j; \; \lambda_j \leq \lambda\}=\frac{\text{Area}\,(\Omega)}{2\pi}\lambda+O(\sqrt\lambda).$$ This asymptotic in particular implies that Area$(\Omega)$ is a spectral invariant. Hilbert thought this conjecture would not be proven in his lifetime,
but less than two years later Hermann Weyl proved it using the theory of integral equations taught to him by Hilbert. In 1954, Pleijel \cite{ple} proved that one knows the perimeter of $\Omega$, and in \cite{Kac}, Kac rephrased these results in terms of asymptotics of the heat trace
\[
\Tr e^{-t\Delta_\Omega} \sim t^{-1}\sum_{j=0}^\infty  a_j t^{j/2}, \qquad t \to 0^+,
\]
where the first coefficient $a_0$ gives the area and the second coefficient gives the perimeter. In 1967, McKean and Singer \cite{ms} proved Pleijel's conjecture that the Euler characteristic $\chi(\Omega)$ is also a spectral invariant (this is in fact given by $a_2$) and hence the number of holes is known. 
In 1991, Gordon, Webb and Wolpert \cite{GWW}, found examples of pairs of distinct plane
domains with the same spectrum. However, their examples were non-convex and non-smooth, and it remains an open question to prove that convex domains are determined by the spectrum (although there are higher dimensional counterexamples for this by Gordon-Webb \cite{GW}) or that smooth domains are determined by the spectrum.

In this survey we review positive inverse spectral and inverse resonant results for the following kinds of problems: Laplacians on bounded domains, Laplace-Beltrami operators on compact manifolds, Schr\"odinger operators, Laplacians on exterior domains, and Laplacians on manifolds which are hyperbolic near infinity. We also recommend the recent survey of Zelditch \cite{zelditchsurvey}. For negative results (counterexamples) we refer the reader to the surveys of Gordon \cite{gordon2000} and Gordon-Perry-Schueth \cite{gordon2005}.

In the next two sections of the paper we review uniqueness results for radial problems \S\ref{s:radial}, and for real analytic and symmetric problems \S\ref{s:analytic}. In the first case the object to be identified satisfies very strong assumptions (radialness includes full symmetry as well as analyticity) but it is identified in a broad class of objects. In this case the first few heat invariants, together with an isoperimetric or isoperimetric-type inequality, often suffice.  In the second case the assumptions on the object to be identified are somewhat weaker (only analyticity and finitely many reflection symmetries are assumed) but the identification is only within a class of objects which also satisfies the same assumptions, and generic nondegeneracy assumptions are also needed. These proofs are based on wave trace invariants corresponding to a single nondegenerate simple periodic orbit and its iterations.

In \S\ref{s:rigidity} we consider rigidity and local uniqueness results, where it is shown in the first case that isospectral deformations of a given object are necessarily trivial, and in the second case that a given object is determined by its spectrum among objects which are nearby in a suitable sense. Here the objects to be determined are more general than in the cases considered in \S\ref{s:radial}, but less general than those in \S\ref{s:analytic}: they are ellipses, spheres, flat manifolds (which have completely integrable dynamics), and manifolds of constant negative curvature (which have chaotic dynamics). The proofs use these special features of the classical dynamics.

In \S\ref{s:compactness} we consider compactness results, where it is shown that certain isospectral families are compact in a suitable topology. These proofs are based on heat trace invariants and on the determinant of the Laplacian, and much more general assumptions are possible than in the previous cases.

Finally, in \S\ref{s:trace} we review the trace invariants used for the positive results in the previous sections, and give examples of their limitations, that is to say examples of objects which have the same trace invariants but which are not isospectral. At this point we also discuss the history of these invariants, going back to the seminal paper of Selberg \cite{selberg}.

We end the introduction by presenting the four basic settings we consider in this survey:

\subsection{Dirichlet and Neumann Laplacians on bounded domains in $\R^n$}
Let $\Omega$ be a bounded open set with piecewise smooth boundary. Let $\Delta_\Omega$ be the nonegative Laplacian on $\Omega$  with Dirichlet or Neumann boundary conditions. Let 
\[
\spec(\Delta_\Omega) = (\lambda_j)_{j=0}^\infty, \qquad  \lambda_0<\lambda_1\le\lambda_2\le \cdots
\]
be the eigenvalues included according to multiplicity, and $u_j$ the corresponding eigenfunctions, that is to say
\[
\Delta_\Omega u_j = \lambda_j u_j.
\]
Recall that $\lambda_0 > 0$ in the Dirichlet case and $\lambda_0 = 0$ in the Neumann case.

\subsection{Laplace-Beltrami operators on compact manifolds}
Let $(M,g)$ be a compact Riemannian manifold without boundary. Let $\Delta_g = - \Div_g \grad_g$ be the nonnegative Laplace-Beltrami operator on $M$, which we also call the Laplacian for short.  Let
\[
\spec(\Delta_g) = (\lambda_j)_{j=0}^\infty, \qquad  0= \lambda_0<\lambda_1\le\lambda_2\le \cdots
\]
 be the eigenvalues included according to multiplicity, and $u_j$ the corresponding eigenfunctions, that is to say
\[
\Delta_g u_j = \lambda_j u_j.
\]

\subsection{Nonsemiclassical and semiclassical Schr\"odinger operators on $\R^n$}
Let 
\begin{equation}\label{e:vcond}
V \in C^\infty(\R^n;\R), \qquad \lim_{|x| \to \infty} V(x) = \infty,
\end{equation}
 and let $\Delta$ be the nonnegative Laplacian on $\R^n$. Let
\[
P_{V,h} = h^2 \Delta + V, \qquad h>0,
\]
\[
P_V = P_{V,1}.
\]
We call $P_V$ the nonsemiclassical Schr\"odinger operator associated to $V$, and $P_{V,h}$ the semiclassical operator. For any $h>0$, the spectrum of $P_{V,h}$ on $\R^n$ is discrete, and we write it as 
\[
\spec(P_{V,h}) = (\lambda_j)_{j=0}^\infty, \qquad  \lambda_0<\lambda_1\le\lambda_2\le \cdots.
\]
The eigenvalues $\lambda_j$ depend on $h$, but we do not include this in the notation. We denote by $u_j$ the corresponding eigenfunctions (which also depend on $h$), so that
\[
P_{V,h} u_j = \lambda_j u_j.
\]

\subsection{Resonance problems for obstacle and potential scattering}

In this subsection we discuss problems where the spectrum consists of a half line of essential spectrum, together with possibly finitely many eigenvalues. In such settings the spectrum contains limited information, but one can often define resonances, which supplement the discrete spectral data and contain more information.

\subsubsection{Obstacle scattering in $\R^n$}
Let $O \subset \R^n$ be a bounded open set with smooth boundary, let $\Omega = \R^n \setminus \overline O$, and suppose that $\Omega$ is connected. Let $\Delta_\Omega$ be the nonnegative Dirichlet or Neumann Laplacian on $\Omega$. Then the spectrum of $\Delta_\Omega$ is continuous and equal to $[0,\infty)$, and so it contains no (further) information about $\Omega$. One way to reformulate the inverse spectral problem in this case is in terms of \textit{resonances}, which are defined as follows. Introduce a new spectral parameter $z = \sqrt\lambda$, with $\sqrt{}$ taken so as to map $\C \setminus [0,\infty)$ to the upper half plane. As $\im z \to 0^+$, $z^2$ approaches $[0,\infty)$ and the resolvent $(\Delta_\Omega - z^2)^{-1}$ has no limit as a map $L^2(\Omega) \to L^2(\Omega)$. However, if we restrict the domain of the resolvent and expand the range it is possible not only to take the limit but also to take a meromorphic continuation to a larger set. More precisely the resolvent
\[
(\Delta_\Omega - z^2)^{-1} \colon  L^2_{\comp} \to L^2_{\loc},
\]
(where $L^2_{\comp}$ denotes compactly supported $L^2$ functions and $L^2_{\loc}$ denotes functions which are locally $L^2$) continues meromorphically as an operator-valued function of $z$ from $\{\im z > 0\}$ to $\C$ when $n$ is odd and to the Riemann surface of  $\log z$ when $n$ is even. Resonances are defined to be the poles of this continuation of the resolvent. Let $\res(\Delta_\Omega)$ denote the set of resonances, included according to multiplicity. See for example \cite{melbook, sjolec, Zworski:resbook} for more information.

\subsubsection{Potential scattering in $\R^n$}
Let $P_{V,h}$ be as before, but instead of \eqref{e:vcond} assume $V \in C_0^\infty(\R^n)$. Then the continuous spectrum of $P_{V,h}$ is equal to $[0,\infty)$, but if $V$ is not everywhere nonnegative then $P_{V,h}$ may have finitely many negative eigenvalues. In either case, the resolvent 
\[
(P_{V,h} - z^2)^{-1} \colon  L^2_{\comp} \to L^2_{\loc},
\]
has a meromorphic continuation from $\{\im z > 0\}$ to $\C$ when $n$ is odd and to the Riemann surface of  $\log z$ when $n$ is even, and resonances are defined to be the poles of this continuation. Let $\res(P_{V,h})$ denote the set of resonances, included according to multiplicity. Again, see for example \cite{melbook, sjolec, Zworski:resbook} for more information.

\subsubsection{Scattering on asymptotically hyperbolic manifolds}

The problem of determining a non-compact manifold from the scattering resonances of the associated Laplace-Beltrami is in general a much more difficult one, but some progress has been made in the asymptotically hyperbolic setting. Meromorphic continuation of the resolvent was established by Mazzeo-Melrose \cite{mm}, and a wave trace formula in the case of surfaces with exact hyperbolic ends was found by Guillop\'e-Zworski \cite{gz}, which has led to some compactness results: see \S\ref{s:manifold}.

\section{The radial case}\label{s:radial}

In this case one makes a strong assumption (radial symmetry) on the object to be spectrally determined (whether it is an open set in $\R^n$, a compact manifold, or a potential) but makes almost no assumption on the class of objects within which it is determined. The methods involved use the first few heat invariants, and in many cases the isoperimetric inequality or an isoperimetric-type inequality.

\subsection{Bounded domains in $\R^n$}\label{s:radialbounded}
The oldest inverse spectral results are for radial problems. If $\Omega \subset \R^n$ is a bounded open set with smooth boundary, then the spectrum of the Dirichlet (or Neumann) Laplacian on $\Omega$ agrees with the spectrum on the unit ball if and only if $\Omega$ is a translation of this ball. This can be proved in many ways; one way is to use heat trace invariants. These are defined to be the coefficients of the asymptotic expansion of the heat trace as $t \to 0^+$:
\begin{equation}\label{e:heattrace}
\sum_{j=0}^\infty e^{-t\lambda_j} = \Tr e^{-t\Delta_\Omega} \sim t^{-n/2}\sum_{j=0}^\infty a_jt^{j/2},
\end{equation}
where in both the Dirichlet and the Neumann case $a_0$ is a universal constant times $\vol(\Omega)$, and $a_1$ is a universal constant times $\vol(\D \Omega)$. The left hand side is clearly determined by the spectrum, and so the conclusion follows from the isoperimetric inequality.

\subsection{Compact manifolds}
Let $(M,g)$ be a smooth Riemannian manifold of dimension $n$ without boundary.  If $n \le 6$, then the spectrum of the Laplacian on $M$ agrees with the spectrum on $S^n$ (equipped with the round metric) if and only if $M$ is isometric to $S^n$. This was proved by Tanno in \cite{t2,t1} using the first four coefficients, $a_0,a_1,a_2,a_3$, of the heat trace expansion, which in this case takes the form
\[
\sum_{j=0}^\infty e^{-t\lambda_j} = \Tr e^{-t\Delta_g} \sim t^{-n/2}\sum_{j=0}^\infty a_jt^{j}.
\]
In higher dimensions the analogous result is not known. In \cite{zel96}, Zelditch proves that if the multiplicities $m_k$ of the \textit{distinct} eigenvalues $0=E_0 < E_1 < E_2< \cdots$ of the Laplacian on $M$ obey the asymptotic $m_k = a k^{n-1} + O(k^{n-2})$, for some $a>0$ as $k \to \infty$ (this is the asymptotic behavior for the multiplicities of the eigenvalues of the sphere), then $(M,g)$ is a Zoll manifold, that is to say a manifold on which all geodesics are periodic with the same period.

\subsection{Schr\"odinger operators}
In general it is impossible to determine a potential $V$ from the spectrum of the nonsemiclassical Schr\"odinger operator $\Delta + V$. For example, in \cite{MT81}, McKean and Trubowitz find an infinite dimensional family of potentials in $C^\infty(\R)$ which are isospectral with the harmonic oscillator $V(x) = x^2$. 

However, analogous uniqueness results to those above are proved in \cite{dhv} by Ventura and the authors, where it is shown that radial, monotonic potentials in $\R^n$ (such as for example the harmonic oscillator) are determined by the spectrum of the associated \textit{semiclassical} Schr\"odinger operator among all potentials with discrete spectrum. The approach is based in part on that of Colin de Verdi\`ere \cite{c}, and Guillemin-Wang \cite{gw10} (see also \cite[\S 10.6]{gs}), where a one dimensional version of the result is proved. They show that an even function (or a suitable noneven function) is determined by its spectrum within the class of functions monotonic away from $0$.

The method of proof is similar to that used to prove spectral uniqueness of balls in $\R^n$ as discussed in \S\ref{s:radialbounded} above. Namely, we use the first two trace invariants, this time of the semiclassical trace formula of Helffer-Robert \cite{hr}, together with the isoperimetric inequality.  We show that if $V,V_0$ are as in \eqref{e:vcond}, if  $V_0(x) = R(|x|)$  where $R(0) = 0$ and $R'(r) > 0$ for $r>0$, and if $\spec(P_{V,h}) = \spec(P_{V_0,h})$ up to order\footnote{The implicit rate of convergence here must be uniform on $[0,\lambda_0]$ for each $\lambda_0>0$.}  $o(h^2)$ for $h \in \{h_j\}_{j=0}^\infty$ with $h_j \to 0^+$, then $V(x) = V_0(x-x_0)$ for some $x_0 \in \R^n$.

The semiclassical trace formula we use is
\begin{equation}\label{e:trace}
\begin{split}\Tr(&f(P_{V,h})) = \\&\frac 1 {(2\pi h)^n} \left(\int_{\R^{2n}} f(|\xi|^2 + V)dxd\xi + \frac{h^2}{12} \int_{\R^{2n}} |\nabla V|^2 f^{(3)}(|\xi|^2 + V)dxd\xi + \Oh(h^4)\right),\end{split}
\end{equation}
where $f \in C_0^\infty(\R)$.

Because the spectrum of $P_{V,h}$ is known up to order $o(h^2)$ we obtain from \eqref{e:trace} the two trace invariants
\begin{equation}\label{e:gsinv}\int_{\{|\xi|^2 + V(x) < \lambda\}} dxd\xi, \qquad \int_{\{|\xi|^2 + V(x) < \lambda\}} |\nabla V(x)|^2dxd\xi,\end{equation}
for each $\lambda$. It follows in particular that $V$ is nonnegative. By integrating in the $\xi$ variable, we rewrite these invariants as follows:
\begin{equation}\label{e:same}\int_{\{V(x) < \lambda\}} (\lambda - V)^{n/2} dx, \qquad \int_{\{V(x) < \lambda\}} |\nabla V(x)|^2 (\lambda - V)^{n/2} dx.\end{equation}

Using the coarea formula we rewrite the invariants in \eqref{e:same} as
\[\int_0^\lambda \left(\int_{\{V=s, \nabla V \ne 0\}} \frac{(\lambda - V)^{n/2}}{|\nabla V|} dS\right)ds, \qquad \int_0^\lambda \left(\int_{\{V=s\}} {|\nabla V|}{(\lambda - V)^{n/2}} dS\right)ds.\]

Using the fact that $V=s$ in the inner integrand, the factor of $(\lambda - V)^{n/2} = (\lambda - s)^{n/2}$ can be taken out of the surface integral, leaving
\begin{equation}\label{e:abel}\int_0^\lambda (\lambda - s)^{n/2} I_1(s)ds, \qquad \int_0^\lambda (\lambda - s)^{n/2} I_2(s)ds,\end{equation}
where
\begin{equation}\label{e:surfaceinvariants}I_1(s) = \int_{\{V=s, \nabla V \ne 0\}} \frac 1 {|\nabla V|} dS, \qquad I_2(s) = \int_{\{V=s\}} |\nabla V| dS. \end{equation}

We denote the integrals \eqref{e:abel} by $A_{1 + n/2}(I_1)(\lambda)$ and $A_{1+ n/2}(I_1)(\lambda)$. These are Abel fractional integrals of $I_1$ and $I_2$ (see for example \cite[\S5.2]{Zelditch} and \cite[(10.45)]{gs}), and they can be inverted by applying $A_{1 + n/2}$, using the formula 
\begin{equation}\label{e:abel2}\frac 1 {\Gamma(\alpha)}A_\alpha\circ \frac 1 {\Gamma(\beta)}A_\beta = \frac 1 {\Gamma(\alpha + \beta)}A_{\alpha+\beta},\end{equation}
and differentiating $n + 1$ times. From this we conclude that the functions $I_1$ and $I_2$ in \eqref{e:surfaceinvariants} are spectral invariants for every $s >0$.

Integrating $I_1$ and using the coarea formula again we find that the volumes of the sets $\{V < s\}$ are spectral invariants:
\begin{equation}\label{e:volumeinvariants}\int_0^s I_1(s')ds' =\int_0^s \int_{\{V=s', \nabla V \ne 0\}} \frac 1 {|\nabla V|} dSds'  = \int_{\{V< s\}} 1 dx.\end{equation}

From Cauchy-Schwarz and the fact that $I_1$ and $I_2$ are spectral invariants we obtain
\begin{equation}\label{e:cs}\left(\int_{\{V=s\}}1 dS\right)^2 \le \int_{\{V=s\}} \frac 1 {|\nabla V|} dS\int_{\{V=s\}} |\nabla V|dS = \int_{\{R=s\}} \frac 1 {R'} dS \int_{\{R=s\}} R'dS,\end{equation}
when $s$ is not a critical value of $V$, i.e. by Sard's theorem for almost every $s \in (0,\lambda_0)$. On the other hand, using the invariants obtained in \eqref{e:volumeinvariants} and the fact that the sets $\{R<s\}$ are balls, by the isoperimetric inequality we find
\begin{equation}\label{e:spheres} \int_{\{R=s\}} 1 dS \le \int_{\{V=s\}}1 dS.\end{equation}
However,
\[\left(\int_{\{R=s\}} 1 dS\right)^2 =  \int_{\{R=s\}} \frac 1 {R'} dS \int_{\{R=s\}} R'dS,\]
because $1/R'$ and $R'$ are constant on $\{R=s\}$. Consequently
\[\int_{\{R=s\}} 1 dS = \int_{\{V=s\}}1 dS,\]
and so $\{V=s\}$ is a sphere for almost every $s$, because only spheres extremize the isoperimetric inequality. Moreover,
\[\left(\int_{\{V=s\}}1 dS\right)^2 = \int_{\{V=s\}} \frac 1 {|\nabla V|} dS \int_{\{V=s\}} |\nabla V|dS,\]
and so $|\nabla V|^{-1}$ and $|\nabla V|$ are proportional on the surface $\{V = s\}$ for almost every $s$, again by Cauchy-Schwarz. Using the equation \eqref{e:cs} to determine the constant of proportionality, we find that
\[|\nabla V|^2 = R'(R^{-1}(s))^2 = (R^{-1})'(s)^{-2} \Def F(s)\]
on $\{V =s\}$. In other words
\begin{equation}\label{e:pde}|\nabla V|^2 = F(V),\end{equation}
for all $x \in V^{-1}(s)$ for almost all $s$. However, because $F(V) \ne 0$ when $V \ne 0$, it follows by continuity that this equation holds for all $x \in V^{-1}((0,\infty))$.

We solve this equation by restricting it to flowlines of $\nabla V$, with initial conditions taken on a fixed level set $\{V = s_0\}$, and conclude that, the level surfaces are not only spheres (as followed from \eqref{e:spheres}) but are moreover spheres with a common center. Hence, up to a translation, $V$ is radial. Since the volumes \eqref{e:volumeinvariants} are spectral invariants, it follows that $V(x) = R(|x|)$.

\subsection{Resonance problems}

We first mention briefly some results for inverse problems for resonances for the nonsemiclassical Schr\"odinger problem when $n=1$. In \cite{Zworski:isopolar} Zworski proves that a compactly supported even potential $V \in L^1(\R)$ is determined from the resonances of $P_V$ among other such potentials, and in \cite{K} Korotyaev shows that a potential which is not necessarily even is determined by some additional scattering data. 

Analogous results to those discussed in \S\ref{s:radialbounded} hold in the case of obstacle scattering. Hassell and Zworski \cite{hz99} show that a ball is determined by its Dirichlet resonances among all compact obstacles in $\R^3$. Christiansen \cite{chr} extends this result to multiple balls, to higher odd dimensions, and to Neumann resonances. As in the other results discussed above, the proofs use two trace invariants and isoperimetric-type inequalities, although the invariants and inequalities are different here. There is also a large literature of inverse scattering results where data other than the resonances are used. A typical datum here is the \textit{scattering phase}: see for example \cite[\S 4.1]{melbook}.

In \cite{dh} we prove the analogue for resonances of the result in the previous section for semiclassical Schr\"odinger operators with discrete spectrum. Let $n \ge 1$ be odd, and let $V_0,V \in C_0^\infty(\R^n;[0,\infty))$. Suppose $V_0(x) = R(|x|)$, and $R'(r)$ vanishes only at $r=0$ and whenever $R(r) = 0$, and suppose that $\res(P_{V_0,h}) = \res(P_{V,h})$, up to order\footnote{The implicit rate of convergence here must be uniform on the disk of radius $\lambda_0$ for each $\lambda_0>0$.} $o(h^2)$, for $h \in \{h_j\}_{j=1}^\infty$ for some sequence $h_j \to 0$. Then there exists $x_0 \in \R^n$ such that $V(x) = V_0(x-x_0)$.

Our proof is, as before, based on recovering and analyzing first two integral invariants of the Helffer-Robert semiclassical trace formula (\cite[Proposition 5.3]{hr}, see also \cite[\S 10.5]{gs}):
\begin{align}\Tr(f(P_{V,h})&) - f(P_{0,h})) = \label{e:hsres}\\\nonumber\frac 1 {(2\pi h)^n} &\left(\int_{\R^{2n}} f(|\xi|^2 + V) - f(|\xi|^2)dxd\xi + \frac{h^2}{12} \int_{\R^{2n}} |\nabla V|^2 f^{(3)}(|\xi|^2 + V)dxd\xi + \Oh(h^4)\right).\end{align}

To express the left hand side of \eqref{e:hsres} in terms of the resonances of $P_{V,h}$, we use Melrose's Poisson formula (\cite{Melrose:Trace}), an extension of the formula of Bardos-Guillot-Ralston (\cite{bgr}):
\begin{equation}\label{e:melrose}
2\Tr\left(\cos(t\sqrt{P_{V,h}}) - \cos(t\sqrt{P_{0,h}})\right) = \sum_{\lambda \in \res(P_{V,h})} e^{-i|t|\lambda}, \qquad t \ne 0,
\end{equation}
where equality is in the sense of distributions on $\R\setminus 0$.

From \eqref{e:melrose}, it follows that if 
\begin{equation}\label{e:moments}\hat g \in C_0^\infty(\R \setminus 0)  \textrm{ is even,}\end{equation}
 then
\begin{equation}\label{e:melrose2}
\Tr(g(\sqrt{-h^2 \Delta + V}) - g(\sqrt{-h^2\Delta})) = \frac 1 {4\pi}\sum_{\lambda \in \res(P_{V,h})} \int_\R e^{-i|t|\lambda}\hat g(t)dt.
\end{equation}

 Now setting the right hand sides of \eqref{e:melrose2} and \eqref{e:hsres} equal and taking $h \to 0$, we find that
\begin{equation}\label{e:traceinv}
\int_{\R^{2n}}  f(|\xi|^2 + V) - f(|\xi|^2)dxd\xi, \qquad  \int_{\R^{2n}} |\nabla V|^2 f^{(3)} (|\xi|^2 + V)dxd\xi
\end{equation}
are resonant invariants (i.e. are determined by knowledge of the resonances up to $o(h^2)$) provided that $f(\tau^2) = g(\tau)$ for all $\tau$ and for some $g$ as in \eqref{e:moments}. Taylor expanding, we write the first invariant as
\[
\sum_{k=1}^m \frac 1 {k!} \int _{\R^n} f^{(k)}(|\xi|^2)d\xi \int_{\R^n} V(x)^k dx + \int_{\R^{2n}}\frac{V(x)^{m+1}} {m!} \int_0^1(1-t)^mf^{(m+1)}(|\xi|^2 + tV(x))dtdxd\xi.
\]
Replacing $f$ by $f_\lambda$, where $f_\lambda(\tau) = f(\tau/\lambda)$ (note that  $g_\lambda(\tau) = f_\lambda(\tau^2)$ satisfies \eqref{e:moments}) gives
\[
\sum_{k=1}^m \lambda^{n/2-k}\frac 1 {k!} \int _{\R^n} f^{(k)}(|\xi|^2)d\xi \int_{\R^n} V(x)^k dx + \Oh(\lambda^{n/2-m-1})
\]
Taking $\lambda \to \infty$ and $m \to \infty$ we obtain the invariants
\[
\int _{\R^n} f^{(k)}(|\xi|^2)d\xi \int_{\R^n} V(x)^k dx,
\]
for every $k \ge 1$.

In \cite[Lemma 2.1] {dh} it is shown that there exists $g$ satisfying \eqref{e:moments} such that if $f(\tau^2) = g(\tau)$, then $ \int _{\R^n} f^{(k)}(|\xi|^2)d\xi \ne 0$, provided $k \ge n$.

This shows that 
\begin{equation}
\label{e:vk}\int_{\R^{n}} V(x)^k dx = \int_{\R^{n}} V_0(x)^k dx
\end{equation}
for every $k \ge n$, and a similar analysis of the second invariant of \eqref{e:traceinv} proves that
\begin{equation}\label{e:vk2}
\int_{\R^{n}} V(x)^k |\nabla V(x)|^2 dx = \int_{\R^{n}} V_0(x)^k |\nabla V_0(x)|^2 dx
\end{equation}
for every $k \ge n$.

We rewrite the invariant \eqref{e:vk} using $V_*dx$, the pushforward of Lebesgue measure by $V$, as
\begin{equation}\label{e:deltasum}
\int_{\R^{n}} V(x)^k dx = \int_\R s^k(V_*dx)_s = i^{k} \widehat{V_*dx}^{(k)}(0).
\end{equation}
Since $V$ and $V_0$ are both bounded functions, the pushforward measures are compactly supported and hence have entire Fourier transforms, and we conclude that
\[
V_*dx = {V_0}_*dx + \sum_{k=0}^{n-1} c_k \delta^{(k)}_0= {V_0}_*dx +  c_0 \delta_0.
\]
For the first equality we used the invariants \eqref{e:deltasum}, and for the second the fact that $V_*dx$ is a measure.
In other words
\[
\vol(\{V > \lambda\}) = \vol(\{V_0 > \lambda\})
\]
whenever $\lambda > 0$. Moreover, this shows that $V_*dx$ is absolutely continuous on $(0,\infty)$, and so by Sard's lemma the critical set of $V$ is Lebesgue-null on $V^{-1}((0,\infty))$. As a result we may  use the coarea formula\footnote{If $n=1$ we put $\int_{\{V=s\}} |\nabla V|^{-1} dS =\sum_{x \in V^{-1}(s)} |V'(x)|^{-1}.$}  to write
\[
V_*dx = \int_{\{V=s\}} |\nabla V|^{-1} dS ds, \qquad \textrm{on }(0,\infty)
\]
and to conclude that 
\[
 \int_{\{V=s\}} |\nabla V|^{-1} dS = \int_{\{V_0=s\}} |\nabla V_0|^{-1} dS
 \]
for almost every $s>0$.
Similarly, rewriting the invariants \eqref{e:vk2} as
\[
\int_{\R^{n}} V(x)^k |\nabla V(x)|^2 dx = \int_\R s^k \int_{\{V=s\}} |\nabla V|dSds, 
\]
we find that
\[
 \int_{\{V=s\}} |\nabla V|dS =  \int_{\{V_0=s\}} |\nabla V_0|dS, \qquad s>0.
\]
From this point on the proof proceeds as in the previous section.

To our knowledge it is not known whether such results hold in even dimensions. The higher dimensional results discussed above all rely on the Poisson formula \eqref{e:melrose} which is only valid for odd dimensions. A similar formula is also true in the obstacle case \cite{Melrose:Polynomial}, although slightly more care is needed in the definition of $\cos(t\sqrt{\Delta_\Omega}) - \cos(t\sqrt{\Delta_{\R^n}})$ because the two operators act on different spaces.
When $n=1$ a stronger trace formula, valid for all $t \in \R$, is known: see for example \cite[page 3]{Zworski:XEDP}. When $n$ is even, because the meromorphic continuation of the resolvent is not to $\C$ but to the Riemann surface of the logarithm,  Poisson formul\ae ~ for resonances are more complicated and contain error terms: see \cite{Sjostrand:Trace, Zworski:Asian}. A proof based on Sj\"ostrand's local trace formula \cite{Sjostrand:Trace} would be of particular interest, firstly because this formula applies in all dimensions and to a very general class of operators, and also because it uses only resonances in a sector (and in certain versions, as in \cite{bony}, resonances in a strip) around the real axis. This would strengthen the known results in odd dimensions as well as proving results in even dimensions, as one would only have to assume that these resonances agreed and not that all resonances do.

\section{The real analytic and symmetric case}\label{s:analytic}

In this case uniqueness results about nonradial objects are obtained, so the assumptions on the object to be determined are weaker. However, the assumptions on the class of objects within which it is determined are much stronger -- in fact they are the same as the assumptions on the object to be determined. The two main assumptions are analyticity and symmetry. In each case wave invariants are used which are microlocalized near certain periodic orbits, as opposed to the non-microlocal heat invariants of the previous section.

\subsection{Bounded domains in $\R^n$}
Here the main tool is the following result of Guillemin-Melrose \cite{gm}. When $\Omega \subset \R^n$ is a bounded, open set with smooth boundary, they prove that $\Tr(\cos(t\sqrt{\Delta_\Omega}))$ is a tempered distribution in $\R$ with the property
\[
\operatorname{sing\, supp} \Tr(\cos(t\sqrt{\Delta_\Omega})) \subset \{0\} \cup \overline{\textrm{Lsp}(\Omega)},
\]
where $\textrm{Lsp}$ denotes the length spectrum, that is to say the lengths of periodic billiard orbits in $\Omega$. Moreover, they show that if $T \in  \textrm{Lsp}(\Omega)$ is of simple length\footnote{This means that only one periodic orbit (up to time reversal), $\gamma_T$, has length $T$.} and  $\gamma_T$ is nondegenerate\footnote{ This means that $\gamma_T$ is transversal to the boundary and $P_T$ the linearized Poincar\'e map of $\gamma_T$, which is the derivative of the first return map, does not have eigenvalue $1$.}, then for $t$ sufficiently near $T$ we have
\begin{equation}\label{e:guilleminmelrose}\begin{split}
 \Tr&\cos(t\sqrt{\Delta_\Omega}) =\\
 & \re\left[ i^{\sigma_T} \frac{T^\sharp}{\sqrt{|\det(I - P_T)|}}(t - T + i0)^{-1} \left(1 + \sum_{j=1}^\infty a_j (t-T)^j\log(t-T + i0)\right)\right] + S(t),
\end{split}\end{equation}
where $S$ is smooth near $T$. Here $T^\sharp$ is the primitive length of $\gamma_T$, which is the length of $\gamma_T$ without retracing, and $\sigma_T$ is the Maslov index of $\gamma_T$ (which can be defined geometrically but which appears here as the signature of the Hessian in the stationary phase expansion of the wave trace). The coefficients $a_j$ are known as \textit{wave invariants}.  Viewing the boundary locally as the graph of a function $f$, they are polynomials in the Taylor coefficients of $f$ at the reflection points of $\gamma_T$. In general there is no explicit formula, but they were computed by Zelditch \cite{Z09} in the special case discussed below.

Assume now $n=2$, with coordinates $(x,y)$. Assume further that 
\begin{enumerate}
\item $\Omega$ is simply connected, symmetric about the $x$-axis, and $\D \Omega$ is analytic on $\{y \ne 0\}$.
\item There is a nondegenerate vertical bouncing ball orbit $\gamma$ of length $T$ such that both $T$ and $2T$ are simple lengths in $\textrm{Lsp}(\Omega)$.
\item The endpoints of $\gamma$ are not critical points of the curvature of $\D \Omega$.
\end{enumerate}
We recall
that a bouncing ball orbit is a 2-link periodic
trajectory of the billiard flow, i.e. a reversible periodic
billiard trajectory that bounces back and forth along a line
segment orthogonal to the boundary at both endpoints. Without loss of generality we may assume that the bouncing ball orbit in assumption (2) is on the $y$-axis.

Zelditch \cite{Z09} proves that if $\Omega$ and $\Omega'$ both satisfy these assumptions, and if $\spec(\Delta_\Omega) = \spec(\Delta_{\Omega'})$ (for either Dirichlet or Neumann boundary conditions), then $\Omega = \Omega'$ up to a reflection about the $y$-axis. This improves a previous result of Zelditch \cite{zelditchprevious} where an additional symmetry assumption is needed, which in turn improves a previous result of Colin de Verdi\`ere \cite{CV} where rigidity is proved in the class of analytic domains with two reflection symmetries. Under the above assumptions, for $\eps > 0$ sufficiently small, there exists a real analytic function $f\colon (-\eps,\eps) \to \R$ such that 
\[
\Omega \cap \{|x| < \eps\} = \{(x,y)\colon |x| < \eps, |y| < f(x)\}.
\]
To prove the theorem it is enough to show that the Taylor coefficients of $f$ at $0$ are determined by $\spec(\Omega)$ (up to possibly replacing $f(x)$ by $f(-x)$). Zelditch does this by writing a formula for the coefficients $a_j$ of \eqref{e:guilleminmelrose}, which are determined by $\spec(\Omega)$, applied to $\gamma$ and to  $\gamma^2$ (the iteration of $\gamma$):
\begin{equation}\label{e:wavecoeffs}\begin{split}
a_j(\gamma^r) = &A_j(r) f^{(2j+2)}(0) + B_j(r) f^{(2j+1)}(0)f^{(3)}(0) \\&+ \left[\textrm{terms containing } f^{(k)}(0) \textrm{ only for } k \le j\right].
\end{split}\end{equation}
Here $A_j(r)$ and $B_j(r)$ are spectral invariants which are determined by the first term of  \eqref{e:guilleminmelrose}. One can show that $(A_j(1),B_j(1))$ as a vector is linearly independent from $(A_j(2),B_j(2))$. Hence, by an inductive argument, if $f^{(3)}(0) \ne 0$, all the coefficients are determined (up to a sign ambiguity for $f^{(3)}(0)$, which corresponds to reflection about the $y$-axis). The condition $f^{(3)}(0) \ne 0$ is equivalent to assumption (3) above, and in \cite[\S6.9]{Z09} Zelditch outlines a possible proof in the case where $f^{(3)}(0) = 0$.

In \cite{HZ10}, Zelditch and the second author prove that bounded analytic
domains $\Omega \subset \R^n$ with $\pm $ reflection symmetries
across all coordinate axes, and  with one axis height fixed  (and
also satisfying some generic non-degeneracy conditions) are
spectrally determined among other such domains. This inverse
result gives a higher dimensional analogue
of the result discussed above from \cite{Z09}, but with $n$ axes of symmetry rather than $n-1$. To
our knowledge, it is the first positive higher dimensional inverse
spectral result for Euclidean domains which is not restricted to
balls. The proof is based as before on \eqref{e:guilleminmelrose} and on formulas for the $a_j(\gamma^r)$, but there are additional algebraic and combinatorial complications coming from the fact that Taylor coefficients must be recovered corresponding to all possible combinations of partial derivatives. These complications are very similar to those that arise for higher dimensional semiclassical Schr\"odinger operators discussed below in \S\ref{s:analyticschrodinger}

\subsection{Compact manifolds}

To our knowledge all uniqueness results in this category are about surfaces of revolution. In \cite{berard} and \cite{gurarie}, B\'erard and Gurarie show that the joint spectrum of $\Delta_g$ and $\D/\D\theta$ (the generator of rotations) of a smooth surface of revolution determines the metric among smooth surfaces of revolution, by reducing the problem to a semiclassical Schr\"odinger operator in one dimension.

In \cite{bruning}, Br\"uning-Heintze show that the spectrum of $\Delta_g$ alone determines the metric of a smooth surface of revolution with an up-down symmetry among such surfaces. They prove that the spectrum of $\Delta_g$ determines the $S^1$-invariant spectrum (but not necessarily the full joint spectrum), allowing them to apply a result of Marchenko \cite[Theorem 2.3.2]{Marchenko} for one-dimensional Schr\"odinger operators. 

In \cite{Zelditch}, Zelditch proves that a convex analytic surface of revolution satisfying a non-degeneracy condition and a simplicity condition is determined uniquely by the spectrum among all such surfaces. He uses analyticity and convexity to show that the spectrum determines the full joint spectrum of $\Delta_g$ and $\D/\D\theta$, reducing the problem to a semiclassical Schr\"odinger operator in one dimension.

\subsection{Schr\"odinger operators}\label{s:analyticschrodinger}

When $n=1$, Marchenko \cite{Marchenko} shows that an even potential is determined by the spectrum of the associated nonsemiclassical Schr\"odinger operator among all even potentials. More specifically, he proves in \cite[Theorem 2.3.2]{Marchenko} that a Schr\"odinger operator on $[0,\infty)$ is determined by knowledge of both the Dirichlet and the Neumann spectrum. The result for even potentials on $\R$ follows from the result on $[0,\infty)$ as follows: Let $V \in C^\infty(\R)$ obey $\lim_{|x| \to \infty} V(x) = \infty$. If $u_j$ is the eigenfunction of  $-\frac{d^2}{dx^2} + V$ corresponding to the eigenvalue $\lambda_j$, then $u_j$ has exactly $j$ zeros and they are all simple (see  \cite[Chapter 2, Theorem 3.5]{bs}). If $V$ is even then every eigenfunction is either odd or even and this result shows that the parity of $u_j$ is the same as the parity of $j$. In particular $(\lambda'_j)_{j=0}^\infty$ with  $\lambda'_j = \lambda_{2j}$ is the spectrum of
\[
-\frac{d^2}{dx^2} + V \textrm{ on } L^2([0,\infty)) \textrm{ with Neumann boundary condition},
\]
and $(\lambda''_j)_{j=0}^\infty$ with  $\lambda''_j = \lambda_{2j+1}$ is the spectrum of
\[
-\frac{d^2}{dx^2} + V \textrm{ on } L^2([0,\infty)) \textrm{ with Dirichlet boundary condition}.
\]
This reduces the problem on $\R$ to the result of Marchenko.

However, noneven potentials may have the same spectrum: indeed McKean-Trubowitz \cite{MT81} construct an infinite-dimensional family of potentials which have the same spectrum as the one-dimensional harmonic oscillator $V(x) = x^2$.

In \cite{gu07}, Guillemin-Uribe  consider potentials $V$  in $\R^n$ which are analytic and even in all variables, which have a unique global minimum $V(0) = 0$, which obey $\lim\inf_{|x| \to \infty} V(x) > 0$, and such that the square roots of the eigenvalues of $\Hess V(0)$ are linearly independent over $\mathbb{Q}$. They show that such potentials are determined by their low lying semiclassical eigenvalues, that is to say by $\spec(P_{V,h}) \cap [0,\eps]$ for any $\eps > 0$. In \cite{h09}, the second author removes the symmetry assumption in the case $n=1$ but assumes $V'''(0) \ne 0$, and for $n \ge 2$ replaces the symmetry assumption by the assumption that $V(x) = f(x_1^2, \dots x_n^2) + x_n^3 g(x_1^2,\dots,x_n^2)$. In \cite{cg} and \cite{c}, Colin de Verdi\`ere-Guillemin and Colin de Verdi\`ere give another proof of this result for the case $n=1$, and in \cite{gu10}, Guillemin-Uribe give another proof in the higher-dimensional case.

The proofs in \cite{gu07,cg,c,gu10} are based on quantum Birkhoff normal forms. These are a quantum version of the Birkhoff normal forms of classical mechanics. In the classical case, one constructs a symplectomorphism which puts a Hamiltonian function into a canonical form in a neighborhood of a periodic orbit. In the quantum case, one constructs a Fourier integral operator associated to this symplectomorphism which puts a pseudodifferential operator which is a quantization of this Hamiltonian into a canonical form, microlocally near the periodic orbit. Quantum Birkhoff normal forms were developed by Sj\"ostrand \cite{sjo} for semiclassical Schr\"odinger operators near a global minimum of the potential. Guillemin \cite{g96} and Zelditch \cite{zel97,zel98} put the Laplace-Beltrami operator on a compact Riemannian manifold into a quantum Birkhoff normal form. Sj\"ostrand-Zworski \cite{szmonodromy} and Iantchenko-Sj\"ostrand-Zworski \cite{ISZ} studied general semiclassical Schr\"odinger operators on a manifold at nondegenerate energy levels.

The proof in \cite{h09} (that the Taylor coefficients of the potential at the bottom of the well are determined by the low-lying eigenvalues) is based on Schr\"odinger trace invariants. These are coefficients of the expansion
\[
\Tr(e^{-itP_{V,h}/h}\chi(P_{V,h})) = \sum_{j=0}^\infty a_j(t)h^j, \qquad h \to 0^+,
\]
where $\chi \in C_0^\infty(\R)$ is $1$ near $0$ and is supported in a sufficiently small neighborhood of $0$. The coefficients $a_j$ in dimension $n=1$ have exactly the form \eqref{e:wavecoeffs} (and in higher dimensions they have the same form as the higher dimensional coefficients of the wave trace on a bounded domain) and hence, once this fact is established, the remainder of the uniqueness proof is the same for both problems.

\subsection{Resonance problems}

The case of an analytic obstacle with two mutually symmetric connected components is treated by Zelditch \cite{z04} by using the singularities of the wave trace generated by the bouncing ball between the two components. In \cite{z07}, Zworski gives a general method for reducing inverse problems for resonances on a noncompact space to corresponding inverse problems for spectra on a compact space.

In \cite{I08} Iantchenko considers potentials $V$  in $\R^n$ which are analytic and even in all variables, which have a unique global maximum at $V(0) = E$, which extend holomorphically to a sector around the real axis and obey $\lim\inf_{|x| \to \infty} V(x) = 0$ in that sector, and such that the square roots of the eigenvalues of $\Hess V(0)$ are linearly independent over $\mathbb{Q}$. He uses the quantum Birkhoff normal form method of \cite{gu07} to recover the Taylor coefficients of the potential at the maximum and to show that potentials $V$ in this class are determined by the resonances in a small neighborhood of $E$.

\section{Rigidity and local uniqueness results}\label{s:rigidity}

In this section we consider results which show nonexistence of nontrivial isospectral deformations.

\subsection{Bounded domains in $\R^n$}\label{s:domainrigid}

In \cite{Marvizi}, Marvizi-Melrose introduce new invariants for strictly convex bounded domains $\Omega \subset \R^2$ based on the length spectrum, associated with the boundary. They show that, for $m \in \N$ fixed,
\[\begin{split}
\sup\{L(\gamma)\colon \gamma \textrm{ is a periodic billiard orbit with } m \textrm{ rotations and } n \textrm{ reflections}\} \\
\sim m L(\D \Omega) + \sum_{k=1}^\infty c_{k,m} n^{-2k}, \qquad n \to \infty, 
\end{split}\]
where $L$ denotes the length. Then they introduce the following \textit{non-coincidence condition} on $\Omega$, which holds for a dense open family (in the $C^\infty$ topology) of strictly convex domains:  suppose there exists $\eps> 0$ such that if $\gamma$ is a closed orbit with $ L(\D \Omega) -\eps < L(\gamma) < L(\D\Omega)$, then $\gamma$ consists of one rotation. They show that under this condition, the coefficients $c_{k,m}$ are spectral invariants, and they use the invariants $c_{1,1}$ and $c_{2,1}$ to construct a two-parameter family of planar domains which are locally spectrally unique (meaning that each domain has a neighborhood in the $C^\infty$ topology within which it is determined by its spectrum). The two parameter family consists of domains which are defined by elliptic integrals, and which resemble ellipses, but which are not ellipses.

In \cite{gm0}, Guillemin-Melrose consider the Laplacian on an ellipse $\Omega$ given by $x^2/a + y^2/b =  1$, with $a>b>0$, and with boundary condition given by 
\begin{equation}\label{e:k}
\D u/ \D n = Ku, \qquad \textrm{on } \D \Omega,
\end{equation}
where $K \in C^\infty(\D \Omega)$ and is even in both $x$ and $y$. They show $K$ is determined by $\spec(\Delta_{\Omega,K})$, where $\Delta_{\Omega,K}$ is the Laplacian on $\Omega$ with boundary condition \eqref{e:k}.

To explain their method, let us introduce some terminology. For $T>0$ the length of a periodic orbit, the fixed point set of $T$, denoted by $Y_T$, is the set of $(q,\eta) \in B^*\D\Omega$, the coball bundle of $\D\Omega$, such that the billiard orbit corresponding to the initial condition $(q,\eta)$ is periodic and has length $T$. For a more general domain there will often be only one periodic orbit of length $T$ (up to time reversal), but  an ellipse, because of the complete integrability of its billiard flow, always has one or several one-parameter families of such orbits.  Guillemin-Melrose prove that, in the case of the ellipse, for any $T$ which is the length of a periodic orbit such that $L(\D\Omega) - T>0$ is  sufficiently small, $Y_T$ has one connected component $\Gamma$ (up to time reversal). This connected component is necessarily a curve which is invariant under the billiard map.  Moreover, they show that the asymptotic expansion of 
\[
\Tr (\cos(t\sqrt{\Delta_{\Omega,K}}) )- \Tr(\cos(t\sqrt{\Delta_{\Omega,0}}))
\]
 in fractional powers of $t-T$ has leading coefficient 
 \begin{equation}\label{e:leray}
 \int_{\Gamma} \frac K {\sqrt{1 - \eta^2}} d \mu_{\Gamma}.
 \end{equation}
  Here $\mu_{\Gamma}$ is the Leray measure on $\Gamma$. They then show that, under the symmetry assumptions, $K$ is determined from a sequence of such integrals for $T_j$ with $T_j$ tending to $L(\D \Omega)$ from below.

In \cite{HZ10}, Zelditch and the second author prove that an ellipse is infinitesimally spectrally rigid among $C^\infty$ domains with the symmetries of the ellipse. This means that if $\Omega_0$ is an ellipse, and  if $\rho_\epsilon$ is a smooth one-parameter family of smooth functions on $\D\Omega_0$ which are even in $x$ and $y$, and if $\Omega_\epsilon$ is a domain whose boundary is defined by
\[
\D \Omega_\epsilon = \{z + \rho_\epsilon(z) n_z\colon  z \in \D \Omega_0\},
\]
and if $\spec(\Omega_0) = \spec(\Omega_\epsilon)$ for $\epsilon \in [0,\epsilon_0)$, then the Taylor expansion of $\rho_\epsilon$ vanishes at $\epsilon = 0$. In particular, if $\rho$ depends on $\epsilon$ analytically, the deformation is constant. The proof uses Hadamard's variational formula for the wave trace:
\begin{equation}\label{e:hadamard}
\left.\frac d {d\epsilon}\right|_{\epsilon = 0} \Tr (\cos(t\sqrt{\Delta_{\Omega_\epsilon}})) =  \frac{t}{2} \int_{\partial \Omega_0} \D_{n_1} \D_{n_2} S_{\Omega_0}(t,z,z)\left(\left.\frac d {d\epsilon}\right|_{\epsilon = 0}\rho_\epsilon(z)\right) dz,
\end{equation}
where $\D_{n_1}$ and $\D_{n_2}$ denote normal derivatives in the first and second variables respectively, $S_{\Omega_0}$ is the kernel of $\sin(t\sqrt{\Delta_{\Omega_0}})/\sqrt{\Delta_{\Omega_0}}$. They then use \eqref{e:hadamard} to prove that for any $T$ in the length spectrum of $\Omega_0$,   the leading
order singularity of the wave trace variation is,
\begin{equation}\label{e:wavevariation}\begin{split}
\left.\frac d {d\epsilon}\right|_{\epsilon = 0} &\Tr (\cos(t\sqrt{\Delta_{\Omega_\epsilon}})) \sim \\
&\frac{t}{2}\; \re \big\{ \big( \sum_{\Gamma \subset Y_T}
C_{\Gamma} \int_{\Gamma}\left(\left.\frac d {d\epsilon}\right|_{\epsilon = 0}\rho_\epsilon\right) \sqrt{1 - |\eta|^2}  d \mu_{\Gamma}
\big)  (t - T+ i 0)^{- \frac{5 }{2} }\big\},
\end{split}\end{equation}
modulo lower order
singularities, where the sum is over the connected components $\Gamma$ of
the set $Y_T$ of periodic points of the billiard map on $B^*\D\Omega_0$ (and its powers) of length $T$, and where $d\mu_{\Gamma}$ is as in \eqref{e:leray}. As  before, if $L(\D \Omega_0) -T>0$ is sufficiently small, there is only one connected component and the sum has only one term. For an isospectral deformation, the left hand side of \eqref{e:wavevariation} vanishes, and hence the integrals
\[
\int_{\Gamma}\left(\left.\frac d {d\epsilon}\right|_{\epsilon = 0}\rho_\epsilon\right) \sqrt{1 - |\eta|^2}  d \mu_{\Gamma}
\]
vanish when $L(\D \Omega_0) -T>0$ sufficiently small. From this point on proceeding as in \cite{gm0} above one can show that $\left.\frac d {d\epsilon}\right|_{\epsilon = 0}\rho_\epsilon = 0$, and reparametrizing the variation one can show that all Taylor coefficients of the variation are $0$. In \cite{HZ10} it is shown that expansions of the form \eqref{e:hadamard} and \eqref{e:wavevariation} hold more generally and in higher dimensions; indeed \eqref{e:hadamard} holds for any $C^1$ variation of any bounded domain, and a version of \eqref{e:wavevariation} holds whenever the fixed point sets $Y_T$ are clean. These formulas may be useful for example in a possible proof of spectral rigidity of ellipsoids.

\subsection{Compact manifolds}

In \cite{t1} Tanno, uses heat trace invariants to show local spectral uniqueness of spheres in all dimensions. This means that there is a $C^\infty$ neighborhood of the round metric on the sphere within which this metric is spectrally determined. In \cite{Kuw}, Kuwabara does this for compact flat manifolds, and in \cite{sharafutdinov}, Sharafutdinov does this for compact manifolds of constant negative curvature.

In \cite{gk1, gk2}, Guillemin-Kazhdan prove that a negatively curved compact manifold $(M,g)$ with simple length spectrum is spectrally rigid if its sectional curvatures satisfy the pinching condition that for every $x \in M$ there is $A(x)>0$ such that $|K/A + 1| < 1/n$, where $K$ is any sectional curvature at $x$ (note that the pinching condition is satisfied for all negatively curved surfaces because in that case there is only one sectional curvature at  each point $x$ and we may take $A(x) = -K$). Spectrally rigid here means if  $g_\epsilon$ is a smooth family of metrics on $M$ with $g_0 = g$ and with $\spec(\Delta_{g_\epsilon}) = \spec(\Delta_{g})$, then $(M,g_\epsilon)$ is isometric to $(M,g)$ for every $\epsilon$. They further use a similar method of proof to establish a spectral uniqueness result for Schr\"odinger operators on these manifolds. The pinching condition was relaxed by Min-Oo in \cite{minoo} and removed by Croke-Sharafutdinov in \cite{cs}, and the result was extended to Anosov surfaces with no focal points by Sharafutdinov-Uhlmann in \cite{shauhl}.

\subsection{Schr\"odinger operators}

In \cite{h11}, the second author considers anisotropic harmonic oscillators: $V(x) = a_1^2 x_1^2 + \cdots + a_n^2 x_n^2$, where the $a_j$ are linearly independent over $\mathbb Q$. It is shown that if $V_\epsilon(x)$ is a smooth deformation of $V(x)$ within the class of $C^\infty$ functions which are even in each $x_j$, and if $\spec(P_{V,1}) = \spec(P_{V_\epsilon,1})$ for all $\epsilon \in [0,\epsilon)$, then the deformation is flat at $\epsilon = 0$, just as in the infinitesimal rigidity result for the ellipse in \S\ref{s:domainrigid}.

\section{Compactness results}\label{s:compactness}

\subsection{Bounded domains in $\R^n$}

In \cite{Melrose:Compact}, Melrose uses heat trace invariants and Sobolev embedding to prove compactness of isospectral sets of domains $\Omega \subset \R^2$ in the sense of the $C^\infty$ topology on the curvature functions in $C^\infty(\D\Omega)$. This result allows the possibilty of a sequence of isospectral domains whose curvatures converge but which `pinch off' in such a way that the limit object is not a domain, but in \cite{mellectures} points out that this possibility can be ruled out using the fact that the singularity of the wave trace at $t=0$ is isolated. 

In \cite{ops}, Osgood-Phillips-Sarnak give another approach to this problem based on the determinant of the Laplacian. This is defined via the analytic continuation of the zeta function
\[
Z(s) = \sum_{j=1}^\infty\lambda_j^{-s}, \qquad \det \Delta_\Omega = e^{-Z'(0)}.
\]
They consider the domain $\Omega$ as the image of the unit disk $D$ under a conformal map $F$, with $e^{2\phi}g_0$ the induced metric on $D$, where $\phi = \log|F'|$ is a harmonic function. Thus $\phi$ is determined by its boundary values, and the topology of \cite{ops} is the $C^\infty$ topology on $\phi|_{\D D}$, and in this case pinching degenerations are ruled out automatically. In \cite{haze99}, Hassell-Zelditch give a nice review of these results and an application of these methods to the compactness problem for isophasal obstacles in $\R^2$.

To our knowledge there is no compactness result in higher dimensions.

\subsection{Compact manifolds}

In \cite{ops88a,ops88b}, Osgood-Phillips-Sarnak extend their determinant methods to the case of surfaces and prove that the set of isospectral metrics on a given Riemannian surface is sequentially compact in the $C^\infty$ topology, up to isometry. In \cite{cy} and in \cite{bpy}, Chang-Yang and Brooks-Perry-Yang  give compactness results for isospectral metrics in a  given conformal class on a three dimensional manifold. In \cite{bpp}, Brooks-Perry-Petersen prove compactness for isospectral families of Riemannian manifolds provided that either the sectional curvatures are all negative, or that there is a uniform lower bound on the Ricci curvatures. In \cite{zhou}, Zhou shows that on a given manifold, the family of isospectral Riemmanian metrics with uniformly bounded curvature is compact, with no restriction on the dimension.

\subsection{Schr\"odinger operators}

In \cite{bruning84}, Br\"uning considers Schr\"odinger operators $\Delta_g + V$ on a compact Riemannian manifold $(M,g)$, where $V \in C^\infty(M)$, and proves that if the dimension $n \le 3$, then any set of isospectral potentials is compact. In higher dimensions he proves the same result under the additional condition that the $H^s$ norm of $V$ for some $s>3(n/2)-2$ is known to be bounded by some constant $C$. In \cite{donnelly}, Donnelly improves this condition to $s>(n/2)-2$, and derives alternative compactness criteria: he shows that isospectral families of nonnegative potentials are compact in dimensions $n \le 9$. If one considers instead $\Delta_g + \gamma V$, he shows that a family of potentials which is isospectral for more than $(n/2) -1$ different values of $\gamma$ is compact. In particular, this implies compactness of families which are isospectral for the semiclassical problem $h^2\Delta_g + V$.

\subsection{Resonance problems}\label{s:manifold}

Let $(X_0,g_0)$ be a conformally compact surface, which is hyperbolic (has constant curvature) outside a given compact set $K_0 \subset X_0$. This means that, if $K_0$ is taken sufficiently large, then $X_0 \setminus K_0$ is a finite disjoint union of funnel ends, which is to say ends of the form
\begin{equation}\label{e:ends}
(0,\infty)_r \times S^1_\theta, \qquad  dr^2 + \ell ^2 \cosh^2(r)d\theta^2,
\end{equation}
 where $\ell \ne 0$ may vary between the funnels. Then the continuous spectrum of $\Delta_{g_0}$ is given by $[1/4,\infty)$, and the point spectrum is either empty or finite and contained in $(0,1/4)$ (and there is no other spectrum). If we introduce the spectral parameter $z = \sqrt{\lambda - 1/4}$, where $\sqrt{}$ is taken to map $\C \setminus [0,\infty)$ to the upper half plane, then the resolvent $(\Delta_g - 1/4 - z^2)^{-1}$ continues meromorphically from $\{\im z > 0\}$ to $\C$ as an operator $L^2_{\textrm{comp}} \to {L^2_\textrm{loc}}$. This meromorphic continuation can be proved by writing a parametrix in terms of the resolvent of the Laplacian on the ends \eqref{e:ends}, which in this case can be written explicitly in terms of special functions: see \cite{mm} for the general construction, and \cite[\S 5]{gzjfa} for a simpler version in this case.

Borthwick-Perry \cite{bp} use a Poisson formula for resonances due to Guillop\'e-Zworski \cite{gz} and a heat trace expansion to show that the set of surfaces which are isoresonant with $(X_0,g_0)$ and for which there is a compact set $K \subset X$ such that $(X_0 \setminus K_0,g_0)$ is isometric to $(X \setminus K, g)$ is compact in the $C^\infty$ topology, improving a previous result of Borthwick-Judge-Perry \cite{bjp}. They also prove related but weaker results in higher dimensions.

\section{Trace invariants and their limitations}\label{s:trace}

\subsection{Bounded domains in $\R^n$}\label{s:tracedomain}

For $\Delta_\Omega$ with $\Omega \subset \R^n$ a bounded smooth domain we have seen two kinds of trace invariants. The first are heat trace invariants, which are the coefficients $a_j$ of the expansion
\[
\Tr e^{-t\Delta_\Omega} \sim t^{-n/2}\sum_{j=0}^\infty  a_j t^{j/2}, \qquad t \to 0^+,
\]
are given by integrals along the boundary of polynomials in the curvature and its derivatives. These are equivalent to the invariants obtained from coefficients of the expansion of the wave trace  $\Tr\cos(t\sqrt{\Delta_\Omega})$ at $t=0$. 

The other kind are wave trace invariants obtained from coefficients of the expansion of the wave trace at the length of a periodic billiard orbit, always assumed to be nondegenerate and usually assumed to be simple. In this case the formula, as already mentioned in \eqref{e:guilleminmelrose}, is
\begin{equation}\label{e:waveexp}\begin{split}
 \Tr&\cos(t\sqrt{\Delta_\Omega}) =\\
 & \re\left[ i^{\sigma_T} \frac{T^\sharp}{\sqrt{\det(I - P_T)}}(t - T + i0)^{-1} \left(1 + \sum_{j=1}^\infty b_j (t-T)^j\log(t-T + i0)\right)\right] + S(t),
\end{split}\end{equation}
where $\gamma_T$ is the simple periodic orbit of length $T$, and where the coefficients $b_j$ are  polynomials in the Taylor coefficients at the reflection points of $\gamma_T$ of the function of which the boundary is a graph.
Because of this requirement on the periodic orbit, positive inverse results of the kind described above, which are based on the wave trace, always require generic assumptions such as nondegeneracy and simple length spectrum. Although there has been some work on the degenerate case, such as \cite{popov} by Popov, it does not seem to have led yet to uniqueness, rigidity, or compactness results. However, in \cite{Marvizi}, Marvizi-Melrose obtain information from invariants at lengths approaching the length of $\D\Omega$ (see \S\ref{s:domainrigid} above for more information).

Another limitation comes from the fact that domains can have the same trace invariants without being isospectral. That is to say, we can construct $\Omega$ and $\Omega'$ such that $\Tr(\cos(t\sqrt{\Delta_\Omega})) - \Tr(\cos(t\sqrt{\Delta_{\Omega'}})) \in C^\infty(\R)$ (recall that the wave trace invariants are the coefficients in the expansion of the wave trace near a singularity, as in \eqref{e:waveexp}), but $\spec(\Delta_\Omega) \ne \spec (\Delta_{\Omega'})$. This was done by Fulling-Kuchment in \cite{fulling}, following a conjecture of Zelditch \cite{zelditchsurvey}, where the following types of domains are considered (these were first introduced by Penrose to study the illumination problem, and then shown by Lifshits to be examples of nonisometric domains with the same length spectrum):

\begin{figure}[htbp]
\includegraphics[width=12cm]{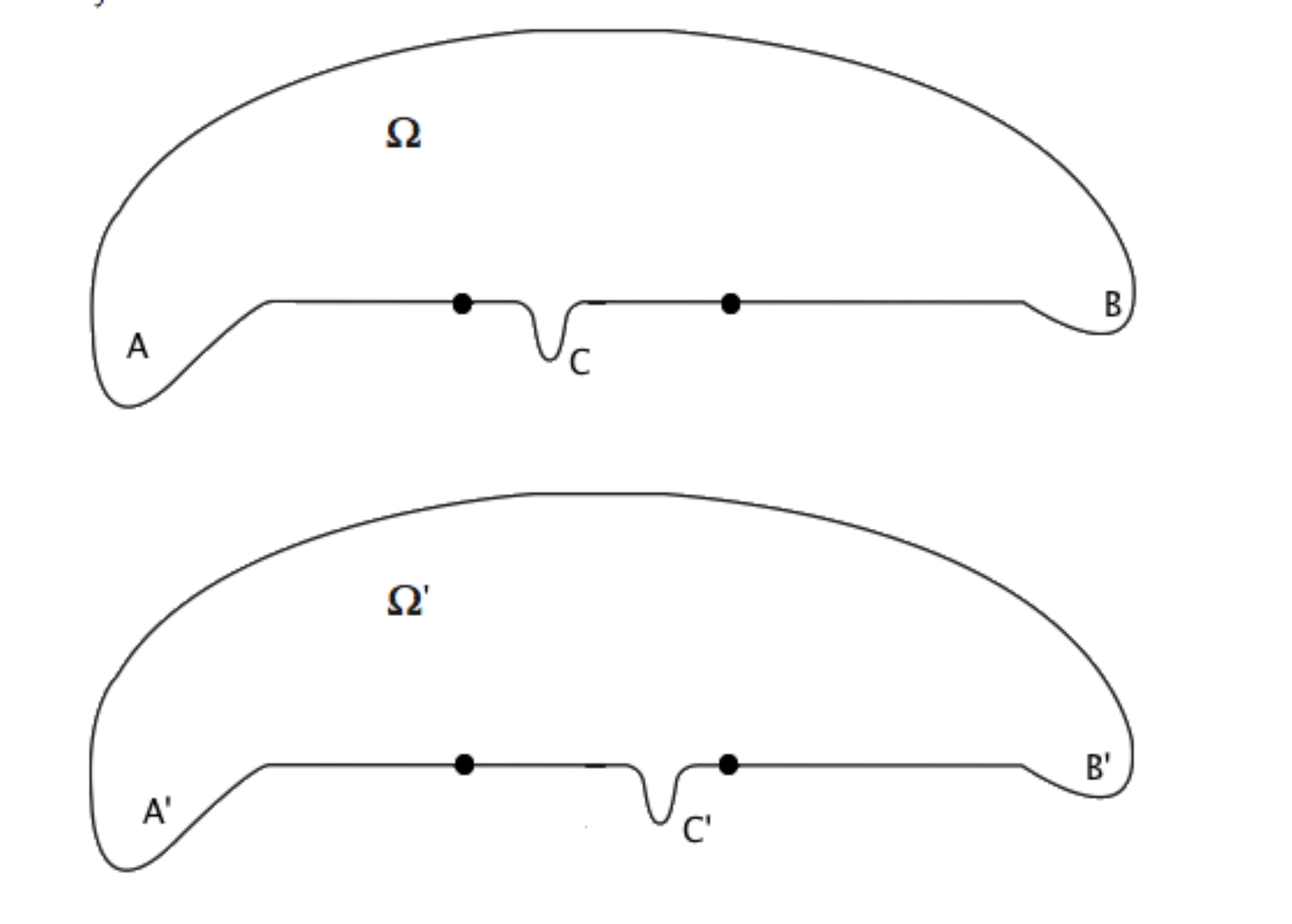}
\caption{Two domains $\Omega$ and $\Omega'$ with $\Tr(\cos(t\sqrt{\Delta_\Omega})) - \Tr(\cos(t\sqrt{\Delta_{\Omega'}})) \in C^\infty(\R)$, but $\spec(\Delta_\Omega) \ne \spec (\Delta_{\Omega'})$. }
\end{figure}

These two domains are obtained by taking a semi-ellipse and adding two asymmetric bumps $A,B$ and $A',B'$, with $A=A'$ and $B = B'$, such that the foci are left unperturbed (as in the figure). Then one adds bumps $C$ and $C'$, the small bumps in the middle which are in between the foci, such that $C \ne C'$ but $C$ and $C'$ are reflections of one another. These two domains are not isometric but have the same heat invariants, because heat invariants are given by integrals along the boundary of polynomials in the curvature and its derivatives -- indeed we have freedom to `slide' $C$ back and forth along the boundary without changing any heat invariants, although this is not the case for the wave trace invariants.

We now show that $\Tr(\cos(t\sqrt{\Delta_\Omega})) - \Tr(\cos(t\sqrt{\Delta_{\Omega'}})) \in C^\infty(\R)$. This is because of the following separation of the phase spaces\footnote{Recall that $B^*\D \Omega$ is the ball bundle, the fibers of which are intervals $[-1,1]$.} $B^*\D \Omega$ and $B^*\D\Omega'$ into two disconnected rooms each, which are invariant under the billiard maps of the domains, and which we denote $R_1$, $R_2$, $R_1'$, and $R_2'$, and which have the property that $\overline{R_1 \cup R_2} = B^*\D\Omega$ and $\overline{R_1' \cup R_2'} = B^*\D\Omega'$. These are defined as follows: $R_1$ is the set of points in $B^*\D\Omega$ whose billiard flowout intersects the part of the boundary strictly in between the two foci, $R_2$ is the set of points in $B^*\D\Omega$ whose billiard flowout intersects the part of the boundary which is strictly outside the two foci but on the axis of the ellipse or below (and similarly for $R_1'$ and $R_2'$).  These two sets are disjoint because billiards in an ellipse which intersect the major axis in between the two foci once do so always. Now we make the generic assumption that no trajectory which passes through the two foci in the initial semi-ellipse is periodic. Because $R_1$ is isometric to $R_1'$, and $R_2$ is isometric to $R_2'$, we have $\Tr(\cos(t\sqrt{\Delta_\Omega})) - \Tr(\cos(t\sqrt{\Delta_{\Omega'}})) \in C^\infty(\R)$. This is because the singularities of $\Tr(\cos(t\sqrt{\Delta_\Omega}))$ occur at $t=T$, where $T$ is the length of a periodic orbit, and only depend on the structure of $B^*\D\Omega$ in an arbitrarily small neighborhood of the orbits of length $T$.
 
To show that $\spec(\Delta_\Omega) \ne \spec(\Delta_{\Omega'})$, Fulling-Kuchment use a perturbation argument based on Hadamard's variational formula for the ground state to show that, for suitably chosen small $C$, the ground states are not the same.

\subsection{Compact manifolds}

As we have already mentioned, heat trace invariants can be defined for compact manifolds $(M,g)$ as well. In the boundaryless case the expansion takes the form
\[
\Tr e^{-t\Delta_g} \sim t^{-n/2}\sum_{j=0}^\infty  a_j t^j, \qquad t \to 0^+,
\]
where the $a_j$ are given by integrals on $M$ of polynomials in the curvature and its derivatives. Half powers of $t$ appear only when there is a boundary, as in the case of domains considered above. Once again, these invariants are equivalent to the invariants obtained from coefficients of the expansion of the wave trace  $\Tr\cos(t\sqrt{\Delta_g})$ at $t=0$. In analogy with the example given in the previous section, we can construct manifolds which are not isometric but which have the same heat invariants by taking a sphere, adding two disjoint bumps, and moving them around. For suitable choices of bumps it should be possible to make the length spectra nonequal, as a result of which the manifolds will be nonisospectral. It seems to be an open problem, however, to find an example of two manifolds $(M,g)$ and $(M',g')$ which are nonisospectral but which have $\Tr(\cos(t\sqrt{\Delta_g})) - \Tr(\cos(t\sqrt{\Delta_{g'}})) \in C^\infty(\R)$, that is to say which have identical wave trace invariants.

In this setting the wave trace expansion was established by Duistermaat-Guillemin \cite{dg}, building off of previous work by Colin de Verdi\`ere \cite{cdvtrace} and Chazarain \cite{chazarain}. It is a generalization of Selberg's Poisson formula \cite{selberg} to an arbitrary compact boundaryless Riemannian manifold. For $T$ the length of a simple nondegenerate periodic geodesic $\gamma_T$, it takes the form
\[
 \Tr e^{it\sqrt{\Delta_g}} =  i^{\sigma_T} \frac{T^\sharp}{\sqrt{|\det(I - P_T)|}}(t - T + i0)^{-1} \left(1 + \sum_{j=1}^\infty b_j (t-T)^j\log(t-T + i0)\right) + S(t),
\]
where $S(t)$ is smooth near $T$.
Using quantum Birkhoff normal forms, in \cite{zel98}, Zelditch shows these coefficients $b_j$ to be integrals of polynomials in the metric and its derivatives along $\gamma_T$. See also \cite{zelwavesurvey} for a more detailed survey on wave invariants. Because of this very local nature of these invariants, to prove uniqueness results one must either assume analyticity (as is done in the results discussed above) or find a way to combine information from many different orbits (no one seems to have been able to do this so far).

\subsection{Schr\"odinger operators}

In the setting of semiclassical Schr\"odinger operators the analogue of the Duistermaat-Guillemin wave trace is the Gutzwiller trace formula near the length $T$ of a periodic trajectory of the Hamiltonian vector field $H_p$ in $p^{-1}(E)$, where $p(x,\xi) = |\xi|^2 + V(x)$:
\[
\Tr e^{-it(P_{V,h}-E)/h} \chi(P_{V,h}) \sim \sum_{\gamma}i^{\sigma_\gamma} \frac{e^{iS_\gamma/h}}{\sqrt{|\det(I - P_\gamma)|}} \sum_{j=0}^\infty a_{j,\gamma}h^j, \qquad a_{0,\gamma} = \delta_0(t-T),
\]
for $t$ near $T$, where $\chi \in C_0^\infty(\R)$ has $\chi = 1$ near $E$. Here the sum in $\gamma$ is over periodic trajectories in $p^{-1}(0)$ of length $T$ and the $a_{j,\gamma}$ are distributions whose singular support is contained in $\{T\}$. This formula  goes back to work of Gutzwiller \cite{gutzwiller}, and was proved in various degrees of generality and with various methods by Guillemin-Uribe \cite{gu89} and Combescure-Ralston-Robert \cite{crr} (see also this last paper for further history and references). This formula is also valid for more general pseudodifferential operators of real principal type, so long as $E$ is a regular value of the principal symbol $p$ and so long as the periodic trajectories in $p^{-1}(E)$ are nondegenerate (so that the determinants in the denominator are nonzero). In particular it also applies on manifolds. In \cite{ISZ}, Iantchenko-Sj\"ostrand-Zworski use it to recover quantum and classical Birkhoff normal forms of semiclassical classical Schr\"odinger operators, at nondegenerate periodic orbits. When the energy level is degenerate, the Gutzwiller trace formula becomes more complicated: see \S\ref{s:analyticschrodinger} for a discussion of the case where $p$ has a unique global minimum at $E$, and see for example \cite{BPU} and \cite{KD} for other cases.

In \cite{c}, Colin de Verdi\`ere gives an example of a pair of potentials $V \not\equiv V' \in C^\infty(\R)$ such that $\spec(P_{V,h}) = \spec(P_{V',h})$ up to $\Oh(h^\infty)$, so that in particular all semiclassical trace invariants for these two potentials agree. He conjectures however, that the spectra are not equal. In \cite{gh11}, Guillemin and the second author construct a pair of potentials, which are perturbations of the harmonic oscillator analogous to the perturbations of the semi-ellipse discussed in \S\ref{s:tracedomain}, and which have different ground states and hence different spectra, although the spectra still agree up to $\Oh(h^\infty)$.


\begin{thebibliography}{00}

\bibitem[BaGuRa82]{bgr} Claude Bardos, Jean-Claude Guillot and James Ralston. La relation de Poisson pour l'\'equation des ondes dans un ouvert non born\'e. Application \`a la th\'eorie de la diffusion. [The Poisson relation for the wave equation in an unbounded open set. Application to scattering theory.] \textit{Comm. Partial Differential Equations} 7:8, 905--958, 1982.

\bibitem[B\'er76]{berard}Pierre B\'erard. Quelques remarques sur les surfaces de r\'evolution dans $\R^3$. [Some remarks on surfaces of revolution in $\R^3$].
\textit{C. R. Acad. Sci. Paris S\'er. A-B} 282:3 (1976) Aii, A159--A161. 

\bibitem[BeSh91]{bs} Felix A. Berezin and Mikhail A. Shubin. The Schr\"odinger equation. \textit{Mathematics and its Applications} 66, 1991.

\bibitem[Bon02]{bony} Jean-Fran\c cois Bony. Minoration du nombre de r\'esonances engendr\'ees par une trajectoire ferm\'ee.  [Lower bound for the number of resonances generated by a closed trajectory]. \textit{Comm. Partial Differential Equations} 27:5, 1021--1078, 2002.

\bibitem[BoJuPe03]{bjp} David Borthwick, Chris Judge, and Peter A. Perry. Determinants of Laplacians and isopolar metrics on surfaces of infinite area. \textit{Duke Math. J.} 118:1, 61--102, 2003.

\bibitem[BoPe11]{bp} David Borthwick and Peter A. Perry. Inverse scattering results for manifolds hyperbolic near infinity. \textit{J. Geom. Anal.} 21:2, 305--333, 2011.

\bibitem[BrPePe92]{bpp} Robert Brooks, Peter Perry, and Peter Petersen. Compactness and finiteness theorems for isospectral manifolds.  \textit{J. Reine Angew. Math.} 426, 67--89, 1992.

\bibitem[BrPeYa89]{bpy} Robert Brooks, Peter Perry, and Paul Yang. Isospectral sets of conformally equivalent metrics. \textit{Duke Math. J.} 58:1, 131--150, 1989.

\bibitem[BrPaUr95] {BPU} Raymond Brummelhuis, Thierry Paul and Alejandro Uribe. Spectral estimates around a critical level. \textit{Duke Math. J.} 78:3, 477--530, 1995.

\bibitem[Br\"u84]{bruning84} Jochen Br\"uning. On the compactness of isospectral potentials. \textit{Comm. Partial Differential Equations} 9:7, 687--698, 1984.

\bibitem[BrHe84]{bruning} Jochen Br\"uning and Ernst Heintze. Spektrale Starrheit gewisser Drehfl\"achen.  [Spectral rigidity of certain surfaces of revolution.] \textit{Math. Ann.} 269:1, 95--101, 1984.

\bibitem[ChYa89]{cy} Sun-Yung A. Chang and Paul C. Yang. Compactness of isospectral conformal metrics on $S^3$. \textit{Comment. Math. Helv.} 64:3, 363--374, 1989.

\bibitem[Cha74]{chazarain} Jacques Chazarain. Formule de Poisson pour les vari\'et\'es riemanniennes. [Poisson formula for Riemannian manifolds]. \textit{Invent. Math.} 24, 65--82, 1974. 

\bibitem[Chr08]{chr} Tanya J. Christiansen. Resonances and balls in obstacle scattering with Neumann boundary conditions. \textit{Inverse Probl. Imaging} 2:3, 335--340, 2008. 

\bibitem[Col73]{cdvtrace} Yves Colin de Verdi\`ere. Spectre du laplacien et longueurs des g\'eod\'esiques p\'eriodiques. II. [Spectrum of the Laplacian and lengths of periodic geodesics]. \textit{Compositio Math.} 27:2, 159--184, 1973.

\bibitem [Col84]{CV} Yves Colin de Verdi\`ere. Sur les longueurs des trajectoires p\'eriodiques d'un billard. [On the lengths of the periodic trajectories of a billiard]. South Rhone seminar on geometry, III (Lyon, 1983), 122--139, \textit{Travaux en Cours}, Hermann, Paris, 1984

\bibitem[Col08p]{c} Yves Colin de Verdi\`ere. A semi-classical inverse problem II: reconstruction of the potential. To appear in \textit{Proc. Duistermaat Conf.} Preprint available at arXiv:0802.1643, 2008.

\bibitem[CoGu08p]{cg} Yves Colin de Verdi\`ere and Victor Guillemin. A semi-classical inverse problem I: Taylor expansions. To appear in \textit{Proc. Duistermaat Conf.}  Preprint available at arXiv:0802.1605, 2008.

\bibitem[CoRaRo99]{crr} Monique Combescure, James Ralston and Didier Robert. A proof of the Gutzwiller semiclassical trace formula using coherent states decomposition. \textit{Comm. Math. Phys.} 202, 463--480, 1999.

\bibitem[CrSh98]{cs} Christopher B. Croke and Vladimir A. Sharafutdinov. Spectral rigidity of a compact negatively curved manifold. \textit{Topology} 37:6, 1265--1273, 1998

\bibitem[DaHe11p]{dh} Kiril Datchev and Hamid Hezari. Resonant uniqueness of radial semiclassical Schr\"odinger operators. To appear in \textit{Appl. Math. Res. Express}. Preprint available at arXiv:1107.0960, 2011.

\bibitem[DaHeVe11]{dhv} Kiril Datchev, Hamid Hezari and Ivan Ventura. Spectral uniqueness of radial semiclassical Schr\"odinger operators. \textit{Math. Res. Lett.} 18:3, 521--529, 2011.


\bibitem[Don05]{donnelly} Harold Donnelly. Compactness of isospectral potentials. \textit{Trans. Amer. Math. Soc.} 357:5, 1717--1730, 2005.

\bibitem[DuGu75]{dg}  Johannes J. Duistermaat and Victor W. Guillemin. The spectrum of positive elliptic operators and periodic bicharacteristics. \textit{Invent. Math.} 29:1, 39--79, 1975.

\bibitem[FuKu05]{fulling} Stephen A. Fulling and Peter Kuchment. Coincidence of length spectra does not imply isospectrality.  \textit{Inverse Problems} 21:4, 1391--1395, 2005.

\bibitem[Gor00]{gordon2000} Carolyn S. Gordon. Survey of isospectral manifolds. \textit{Handbook of differential geometry} 1, 747--778,  2000.

\bibitem[GoPeSc05]{gordon2005} Carolyn S. Gordon, Peter Perry and Dorothee Schueth. Isospectral and isoscattering manifolds: a survey of techniques and examples. In Geometry, spectral theory, groups, and dynamics. \textit{Contemp. Math.} 387, 157--179, 2005.

\bibitem[GoWe94]{GW} Carolyn S. Gordon and  David L. Webb. Isospectral convex domains in Euclidean space. \textit{Math. Res. Lett.} 1:5, 539--545, 1994.

\bibitem[GoWeWo92]{GWW} Carolyn Gordon, David Webb and Scott Wolpert. Isospectral plane domains and surfaces via Riemannian orbifolds. \textit{Invent. Math.} 110:1, 1--22, 1992.

\bibitem[Gui96]{g96} Victor Guillemin. Wave-trace invariants. \textit{Duke Math. J.} 83:2, 287--352, 1996.

\bibitem[GuHe11p]{gh11} Victor Guillemin and Hamid Hezari. A Fulling-Kuchment theorem for the 1D harmonic oscillator. Preprint, 2011.

\bibitem[GuKa80a]{gk1}Victor Guillemin and David Kazhdan. Some inverse spectral results for negatively curved 2-manifolds. \textit{Topology} 19:3, 301--312, 1980.

\bibitem[GuKa80b]{gk2} Victor Guillemin and David Kazhdan. Some inverse spectral results for negatively curved $n$-manifolds, in Geometry of the Laplace operator.  \textit{Proc. Sympos. Pure Math.}, 36, 153--180, 1980.

\bibitem[GuMe79a]{gm0} Victor Guillemin and Richard Melrose. An inverse spectral result for elliptical regions in $\R^2$. \textit{Adv. in Math.} 32:2, 128--148, 1979.

\bibitem[GuMe79b]{gm} Victor Guillemin and Richard Melrose. The Poisson summation formula for manifolds with boundary. \textit{Adv. in Math.} 32:3 204--232, 1979. 

\bibitem[GuPa10]{gp} Victor Guillemin and Thierry Paul. Some remarks about semiclassical trace invariants and quantum normal forms. \textit{Comm. Math. Phys.} 294:1, 1--19, 2010. 

\bibitem[GuSt11p]{gs} Victor Guillemin and Shlomo Sternberg. Semi-classical analysis. Lecture notes available online at \url{http://math.mit.edu/~vwg/semiclassGuilleminSternberg.pdf}

\bibitem[GuUr89]{gu89} Victor Guillemin and Alejandro Uribe. Circular symmetry and the trace formula. \textit{Invent. Math.} 96, 386--423, 1989.

\bibitem[GuUr07]{gu07} Victor Guillemin and Alejandro Uribe. Some inverse spectral results for semi-classical Schr\"o-dinger operators. \textit{Math. Res. Lett.} 14:4, 623--632, 2007. 

\bibitem[GuUr11]{gu10} Victor Guillemin and Alejandro Uribe. Some inverse spectral results for the two-dimensional Schr\"odinger operator. In Geometry and Analysis.  \textit{Adv. Lect.  Math. (ALM)}, 17:1, 319--328, 2011.

\bibitem[GuWa09p]{gw10} Victor Guillemin and Zuoqin Wang. Semiclassical spectral invariants for Schr\"odinger operators. Preprint available at arXiv:0905.0919, 2009.

\bibitem[GuZw95]{gzjfa} Laurent Guillop\'e and Maciej Zworski. Upper bounds on the number of resonances for non-compact Riemann surfaces. \textit{J. Funct. Anal.} 129:2, 364--389, 1995.

\bibitem[GuZw97]{gz} Laurent Guillop\'e and Maciej Zworski. Scattering asymptotics for Riemann surfaces. \textit{Ann. of Math.(2)} 145:3, 597--660, 1997.

\bibitem[Gur95]{gurarie} David Gurarie. Semiclassical eigenvalues and shape problems on surfaces of revolution. \textit{J. Math. Phys.} 36:4, 19--34, 1995.

\bibitem[Gut71]{gutzwiller} Martin Gutzwiller. Periodic orbits and classical quantization conditions. \textit{J. Math. Phys.} 12, 343--358, 1971.

\bibitem[HaZe99]{haze99} Andrew Hassell and Steve Zelditch. Determinants of Laplacians in exterior domains. \textit{Internat. Math. Res. Notices} 1999:18, 971--1004, 1999.

\bibitem[HaZw99]{hz99} Andrew Hassell and Maciej Zworski. Resonant rigidity of $S^2$. \textit{J. Funct. Anal.} 169:2, 604--609, 1999.

\bibitem[HeRo83]{hr} Bernard Helffer and Didier Robert. Calcul fonctionnel par la transformation de Mellin et op\'erateurs admissibles. [Functional calculus by the Mellin transform and admissible operators]. \textit{J. Funct. Anal.} 53:3, 246--268, 1983. 


\bibitem[Hez09]{h09} Hamid Hezari. Inverse spectral problems for Schr\"odinger operators. \textit{Comm. Math. Phys.} 288:3, 1061--1088, 2009.

\bibitem[Hez11p]{h11} Hamid Hezari. Spectral rigidity of the $n$-dimensional harmonic oscillator. In preparation.

\bibitem[HeZe10]{HZ10a} Hamid Hezari and Steve Zelditch. Inverse spectral problem for analytic $(\Z/2\Z)^n$-symmetric domains in $\R^n$. \textit{Geom. Funct. Anal.} 20:1, 160--191, 2010.

\bibitem[HeZe10p]{HZ10} Hamid Hezari and Steve Zelditch. $C^\infty$ spectral rigidity of the ellipse. Preprint available at arXiv:1007.1741, 2010.

\bibitem[Ian08]{I08} Alexei Iantchenko. An inverse problem for trapping point resonances. \textit{Lett. Math. Phys.} 86:2-3, 151--157, 2008.

\bibitem[IaSjZw02]{ISZ} Alexei Iantchenko, Johannes Sj\"ostrand and Maciej Zworski. Birkhoff normal forms in semi-classical inverse problems. \textit{Math. Res. Lett.} 9:2, 337--362, 2002.

\bibitem[Kac66]{Kac} Mark Kac. Can one hear the shape of a drum? \textit{Amer. Math. Monthly} 73:4 part II,
1--23, 1966.

\bibitem[Khu97] {KD} David Khuat-Duy. A semi-classical trace formula for Schr\"odinger operators in the case of a critical energy level. \textit{J. Funct. Anal.} 146:2, 299--351, 1997.

\bibitem[Kor05]{K} Evgeny Korotyaev. Inverse resonance scattering on the real line. \textit{Inverse Problems} 21:1, 325--341, 2005.

\bibitem[Kuw80]{Kuw} Ruishi Kuwabara. On the characterization of flat metrics by the spectrum. \textit{Comment. Math. Helv.} 55:3, 427--444, 1980.

\bibitem[Mar52]{Marchenko}
\foreignlanguage{russian}{\textrm{Vladimir A. Marchenko. Nekotorye voprosy teorii odnomernyh line\U{i}nyh differentsialp1-nyh operatorov vtorogo poryadka.} 1. \textit{Tr. Mosk. Mat. Ob.},} 
[Vladimir A. Marchenko. Some questions of the theory of one-dimensional linear differential operators of the second order. 1.  \textit{Tr. Mosk. Mat. Ob.}], 1, 327--420, 1952. 

\bibitem[MaMe82]{Marvizi} Shahla Marvizi and Richard Melrose. Spectral invariants of convex planar regions. \textit{J. Differential Geom.} 17:2, 475--502, 1982.

\bibitem[MaMe87]{mm} Rafe R. Mazzeo and Richard Melrose. Meromorphic extension of the resolvent on complete spaces with asymptotically constant negative curvature. \textit{J. Funct. Anal.} 75:2, 260--310, 1987.

\bibitem[McSi67]{ms} Henry P. McKean Jr. and Isadore M. Singer. Curvature and the eigenvalues of the Laplacian. \textit{J. Differential Geometry} 1:1, 43--69, 1967.

\bibitem[McTr81]{MT81} Henry P. McKean and Eugene Trubowitz. The spectral class of the quantum-mechanical harmonic oscillator. \textit{Comm. Math. Phys.} 82:4, 471--495, 1981.

\bibitem[Mel82]{Melrose:Trace} Richard Melrose. Scattering theory and the trace of the wave group. \textit{J. Func. Anal.} 45:1, 29--40, 1982.

\bibitem[Mel83p]{Melrose:Compact} Richard Melrose. Isospectral sets of drumheads are compact in $C^\infty$. Available online at \url{http://www-math.mit.edu/~rbm/papers/isospectral/isospectral.pdf}, 1983.

\bibitem[Mel83]{Melrose:Polynomial} Richard Melrose. Polynomial bound on the number of scattering poles. \textit{J. Func. Anal.} 53:3, 287--303, 1983.

\bibitem[Mel95]{melbook} Richard Melrose. Geometric scattering theory. \textit{Cambridge University Press}, 1995.

\bibitem[Mel96]{mellectures} Richard Melrose. The inverse spectral problem for planar domains, in Instructional Workshop on Analysis and Geometry, Part I. \textit{Proc. Centre Math. Appl. Austral. Nat. Univ.} 34, 137--160, 1996.

\bibitem[Min86]{minoo} Maung Min-Oo. Spectral rigidity for manifolds with negative curvature operator, in Nonlinear problems in geometry, \textit{Contemp. Math.} 51, 99--103, 1986. 

\bibitem[OsPhSa88a]{ops88a}  Brad Osgood, Ralph Phillips and Peter Sarnak. Extremals of determinants of Laplacians. \textit{J. Funct. Anal.} 80:1, 148--211, 1988.

\bibitem[OsPhSa88b]{ops88b}  Brad Osgood, Ralph Phillips and Peter Sarnak. Compact isospectral sets of surfaces. \textit{J. Funct. Anal.} 80:1, 212--234, 1988.

\bibitem[OsPhSa89]{ops} Brad Osgood, Ralph Phillips and Peter Sarnak. Moduli space, heights and isospectral sets of plane domains. \textit{Ann. of Math. (2)} 129:2, 293--362, 1989.


\bibitem[PeWoRe07] {PWR} Niklas Peinecke, Franz-Erich Wolter and Martin Reuter. Laplace spectra as fingerprints for image recognition. \textit{Computer-Aided Design} 39:6, 460--476, 2007.

\bibitem[Ple54]{ple} \AA ke Pleijel. A study of certain Green's functions with applications in the theory of vibrating membranes. \textit{Ark. Mat.} 2, 553--569, 1954.

\bibitem[Pop98]{popov} Georgi Popov. On the contribution of degenerate periodic trajectories to the wave-trace. \textit{Comm. Math. Phys.} 196:2, 363--383, 1998.

\bibitem[Reu07] {R} Martin Reuter. Can one hear Shape? \textit{PAMM Proceedings of GAMM07 and ICIAM07} 7:1 1011101--1011102, 2007.
    
\bibitem[ReWoPe07] {RWP} Martin Reuter, Franz-Erich Wolter and Niklas Peinecke.  Laplace-Beltrami spectra as ``Shape-DNA'' of surfaces and solids. \textit{Computer-Aided Design} 38:4, 342--366, 2007.

\bibitem[RWSN09] {RWShN} Martin Reuter, Franz-Erich Wolter, Martha Shenton and Marc Niethammer. Laplace-Beltrami eigenvalues and topological features of eigenfunctions for statistical shape analysis. \textit{Computer-Aided Design} 41:10, 739--755, 2009.



\bibitem[Sel56]{selberg} Atle Selberg. Harmonic analysis and discontinuous groups in weakly symmetric Riemannian spaces with applications to Dirichlet series. \textit{J. Indian Math. Soc. (N.S.)} 20, 47--87, 1956.

\bibitem[Sha09]{sharafutdinov} Vladimir A. Sharafutdinov. Local audibility of a hyperbolic metric. \textit{Sibirsk. Mat. Zh.} 50:5, 1176--1194, 2009; translation in \textit{Sib. Math. J.} 50:5, 929--944, 2009.

\bibitem[ShUh00]{shauhl} Vladimir Sharafutdinov and Gunther Uhlmann. On deformation boundary rigidity and spectral rigidity of Riemannian surfaces with no focal points. \textit{J. Differential Geom.} 56:1, 93--110, 2000.

\bibitem[Sj\"o92]{sjo} Johannes Sj\"ostrand. Semi-excited states in nondegenerate potential wells. \textit{Asymptotic Anal.} 6:1, 29--43, 1992.

\bibitem[Sj\"o97]{Sjostrand:Trace} Johannes Sj\"ostrand. A trace formula and review of some estimates for
resonances, in Microlocal analysis and spectral theory. \textit{NATO Adv. Sci. Inst. Ser. C Math. Phys. Sci.} 490, 377--437, 1997.

\bibitem[Sj\"o02]{sjolec} Johannes Sj\"ostrand. Lectures on resonances. Available online at \url{http://www.math.polytechnique.fr/~sjoestrand/CoursgbgWeb.pdf}. 

\bibitem[SjZw02]{szmonodromy} Johannes Sj\"ostrand and Maciej Zworski. Quantum monodromy and semi-classical trace formul\ae. \textit{J. Math. Pures Appl.} 81:1, 1--33, 2002.

\bibitem[Tan73]{t2} Shukichi Tanno. Eigenvalues of the Laplacian of Riemannian manifolds. \textit{Tohoku Math. J. (2)} 25, 391--403, 1973.

\bibitem[Tan80]{t1} Shukichi Tanno. A characterization of the canonical spheres by the spectrum. \textit{Math. Z.} 175:3, 267--274, 1980.


\bibitem[Zel96]{zel96} Steven Zelditch. Maximally degenerate Laplacians. \textit{Ann. Inst. Fourier (Grenoble)} 46:2, 547--587, 1996.

\bibitem[Zel97]{zel97} Steve Zelditch. Wave invariants at elliptic closed geodesics. \textit{Geom. Funct. Anal.} 7:1, 145--213, 1997.

\bibitem[Zel98a]{zel98} Steve Zelditch. Wave invariants for non-degenerate closed geodesics. \textit{Geom. Funct. Anal.} 8:1, 179--217, 1998.

\bibitem[Zel98b]{Zelditch} Steve Zelditch. The inverse spectral problem for surfaces of revolution. \textit{J. Differential Geom.} 49:2, 207--264, 1998.

\bibitem[Zel99]{zelwavesurvey} Steve Zelditch. Lectures on wave invariants, in Spectral theory and geometry.
\textit{London Math. Soc. Lecture Note Ser.}, 273, 284--328, 1999. 

\bibitem[Zel00]{zelditchprevious} Steve Zelditch.  Spectral determination of analytic bi-axisymmetric plane domains. \textit{Geom. Funct. Anal.} 10:3, 628--677, 2000.

\bibitem[Zel04a]{z04} Steve Zelditch. Inverse resonance problem for $\Z_2$-symmetric analytic obstacles in the plane. Geometric methods in inverse problems and PDE control. \textit{IMA Vol. Math. Appl.} 137,  289--321, 2004. 

\bibitem[Zel04b]{zelditchsurvey} Steve Zelditch, with an appendix by Johannes Sj\"ostrand and Maciej Zworski. The inverse spectral problem.  \textit{Surv. Differ. Geom.} 9, 401--467, 2004. 

\bibitem[Zel09]{Z09} Steve Zelditch. Inverse spectral problem for analytic domains. II. $\Z_2$-symmetric domains. \textit{Ann. of Math. (2)} 170:1, 205--269, 2009. 

\bibitem[Zho97]{zhou} Gengqiang Zhou. Compactness of isospectral compact manifolds with bounded curvatures. \textit{Pacific J. Math.} 181:1, 187--200, 1997.

\bibitem[Zwo96]{Zworski:XEDP} Maciej Zworski. Poisson {formul\ae} for resonances. \textit{S\'emin. \'Equ. D\'eriv. Partielles} Expos\'e $\textrm{n}^\textrm{o}$ XIII, 12p.,  1996--1997.

\bibitem[Zwo98]{Zworski:Asian} Maciej Zworski. Poisson formula for resonances in even dimensions. \textit{Asian J. Math.} 2:3,  609--617, 1998.

\bibitem[Zwo01]{Zworski:isopolar} Maciej Zworski. A remark on isopolar potentials. \textit{SIAM J. Math. Anal.} 32:6, 1324--1326, 2001.

\bibitem[Zwo07]{z07} Maciej Zworski. A remark on: ``Inverse resonance problem for $\Z_2$-symmetric analytic obstacles in the plane'' by S. Zelditch. \textit{Inverse Probl. Imaging} 1:1 225--227, 2007. 

\bibitem[Zwo11]{Zworski:resbook} Maciej Zworski. Lectures on scattering resonances. Lecture notes available onlline at \url{http://math.berkeley.edu/~zworski/res.pdf}.
\end{thebibliography}
\end{document}